\documentclass[11pt,a4paper]{article}
\usepackage{hyperref}
\usepackage{amsmath,amsfonts,amssymb,theorem,graphicx,verbatim,multirow,fullpage}
\usepackage{colortbl}
\usepackage{booktabs}
\usepackage{longtable}
\usepackage{pb-diagram}
\usepackage{algorithm}

\restylefloat{table}
\newcommand{\CC}{\mathbb C}

\newcommand{\AAA}{\mathbb A}
\newcommand{\PP}{\mathbb P}
\newcommand{\ZZ}{\mathbb Z}

\newcommand{\QQ}{\mathbb Q}

\newcommand{\lam}{\lambda}

\newcommand{\Lam}{\Lambda}
\newcommand{\al}{\alpha}
\newcommand{\be}{\beta}
\newcommand{\gog}{\mathfrak g}

\newcommand{\la}{\lambda}
\newcommand{\om}{\omega}

\newcommand{\ep}{\varepsilon}
\newcommand{\Si}{\Sigma}

\newcommand{\Oh}{\mathcal O}

\newcommand{\Cc}{\mathcal C}
\newcommand{\sB}{\mathcal B}

\newcommand{\into}{\hookrightarrow}

\theorembodyfont{\rmfamily}
 
 \newtheorem{theorem}[subsection]{Theorem}

 \newtheorem{coro}[subsection]{Corollary}

\newtheorem{example}[subsection]{Example}
\newtheorem{rmk}[subsection]{Remark}
{
\theorembodyfont{\rmfamily}

 \newtheorem{nothing}[subsection]{}
}
 
 \newenvironment{proof}{\paragraph{Proof}}{\par\medskip}

\theorembodyfont{\rmfamily}

\DeclareMathOperator{\Gr}{Gr }

\DeclareMathOperator{\Hom}{Hom}

\DeclareMathOperator{\Proj}{Proj}
\DeclareMathOperator{\Spec}{Spec}

\newtheorem{thm}{Theorem}[section]

\numberwithin{equation}{section}
\usepackage[normalem]{ulem}
\begin{document}

\title{ Computing isolated orbifolds  in weighted flag varieties }
\date{}
\author{ Muhammad Imran Qureshi}

\maketitle
\begin{abstract}
Given a  weighted flag variety \(w\Si(\mu,u)\) corresponding to  chosen fixed parameters \(\mu\) and \(u\), we present an algorithm to compute lists of   all possible   projectively Gorenstein   \(n\)-folds,  having canonical weight \(k\)  and   isolated orbifold points, appearing   as  weighted complete intersections in \(w\Si(\mu,u) \) or  some projective cone(s) over \(w\Si(\mu,u)\).  We apply our algorithm  to compute  lists of interesting classes  of polarized  3-folds with isolated orbifold points in  the  codimension 8 weighted \(G_2\) variety.  We also show the existence of  some  families of log-terminal     \(\QQ\)-Fano 3-folds  in codimension 8 by explicitly  constructing them as quasilinear sections of a weighted \(G_2\)-variety.   
\end{abstract}

\section{Introduction} 
This article has two main parts. In the first part, for a weighted flag variety \(w\Si(\mu,u)\)  corresponding to the  fixed chosen parameters \(\mu \text{ and }u\), we give  an algorithm to compute  all possible families of   isolated orbifolds of fixed dimension \(n\) and canonical weight \(k\)  that are complete intersections in \(w\Si(\mu,u)\) or  some projective cone(s) over \(w\Si(\mu,u)\). The dimension \(n\) of isolated orbifolds is required to be greater than or equal to   \( 2\), due to Theorem  \ref{hs:brz}, and \(k\) is an integer such that the canonical divisor class \(K_X\backsim kD\), for an ample \(\QQ\)-Carter divisor \(D\).  In the second part, we use the algorithm to  compute numerical candidate lists of different types of families of    isolated 3-dimensional orbifolds whose general member  is a quasilinear section of a certain weighted flag variety \(w\Si(\mu,u)\), also called format. In particular, we obtain lists of canonical 3-folds, Calabi--Yau 3-folds and log-terminal \(\QQ\)-Fano 3-folds. The lists are obtained by   implementing  the algorithm in the computer algebra system {\sc Magma} \cite{magma}.  

We are interested in computing lists of well-formed and quasi-smooth projectively Gorenstein  polarized  \(n\)-folds  $(X,D)$, where $X$ is polarized by an ample \(\QQ\)-Cartier  divisor \(D\), appearing as quasilinear sections of  some fixed weighted flag variety \(w \Si(\mu,u)\). The algorithm to compute such lists  of candidate families of \(n\)-folds is  primarily based on   the Hilbert series formula of Buckley--Reid--Zhou \cite[Thrm. 1.3]{brz}. The formula gives the decomposition of the Hilbert series  \(P_X(t)=P_I(t)+P_Q(t)\): \(P_I(t)\) represents the smooth part and \(P_Q(t)\) represents the orbifold part of the \(n\)-fold \((X,D)\).  The term \(P_{Q}(t) \) will have the form  \(\sum_{Q_i\in \sB} P_{Q_i}(t),\)  where \(\sB\) is a collection of singularities (possibly with repeats) known as basket, coming from the embedding of  \(X\) in some weighted projective space. 

The algorithm starts with  computing the Hilbert series and  the canonical divisor class \(\Oh_{w\Si}(-p)\) of some weighted flag variety \(w\Si(\mu,u)\) of dimension \(d\), where  \(p>>0\) is a positive integer. The flag varieties are known to be Fano varieties, whence  \(-p\) for the canonical weight of \(w\Si\).  Then we find all possible \(n\)-fold quasilinear sections \(X\) of \(w\Si\) or of projective cone(s) over \(w\Si\) such that \(\Oh_X(K_X)=\Oh_X(kD)\). For a chosen \(n\)-fold model  \(X\), we find  its Hilbert series  \(P_{X}(t)\) and the corresponding initial term \(P_I(t)\). Since we are interested in quasi-smooth orbifolds, the singularities of \(X\) come from the singularities of the weighted projective space containing \(X\).  We compute  all the possible baskets \(\sB\) of isolated singularities and their contributions \(P_{Q_{i}}(t)\) to the Hilbert series \(P_X(t)\). In the last step  we run through each basket to determine the multiplicities of the orbifold terms \(P_{Q_i}(t)\): the values of \(m_{i}\) in the equation \begin{equation}P_X(t)=P_I(t)+\sum m_iP_{Q_i}(t).\end{equation}     
 Then \(X\) is a suitable candidate for an isolated polarized   \(n\)-fold if \(m_i\ge 0\) for all \(i\).
In fact, the algorithm can be used to find lists of isolated \(n\)-folds  which may be realized as weighted complete intersection in any ambient  weighted projective variety having computable canonical divisor class and Hilbert series. 

We implement  the algorithm in the computer algebra system  \textsc{Magma}  to   compute  lists of candidate families of isolated 3-folds in two cases: the codimension 8 weighted \(G_{2}\) variety and the codimension 3 weighted \(\Gr(2,5)\).  We compute the numerical candidate families of log-terminal \(\QQ\)-Fano 3-folds, Calabi--Yau 3-folds and canonical 3-folds whose general member is a weighted complete intersection in the corresponding ambient variety. The lists are computed by ordering the input parameters \((\mu,u)\)  in order of increasing the sum of the weights \(W\) of the embedding containing \(w\Si.\) The lists do not present  the full classification of such 3-folds in a chosen weighted flag variety \(w\Si\) but are complete up to a certain value of \(W\) in each case.  The list of   33 candidate log-terminal \(\QQ\)-Fano 3-folds has been explicitly checked by using their defining equations;  6 of them exist as actual polarized 3-folds with the desired properties, given in Theorem \ref{thrm:g2fano}. The list also confirms the non-existence of terminal \(\QQ\)-Fano 3-folds in the  weighted \(G_2\) variety. The  lists of families of Calabi--Yau 3-folds and canonical 3-folds  shall only be considered as numerical candidate 3-folds in the \(G_2\) case. Their  existence can be checked by using their   defining equations under their
Pl\"ucker-style embeddings, and checking their type of singularities using those  equations.  The case \(w\Gr(2,5)\) has been included to check the algorithm against already existing lists of  polarized 3-folds.

The main motivation to constructing polarized varieties as  complete intersections in weighted flag varieties \(w\Si\) comes from Mukai's linear section theorem \cite{mukai,mukai4}: every prime Gorenstein Fano 3-fold of genus \(7\le g\le 10\) is a  linear section of some  flag variety. The idea was first generalized to the weighted case by Corti--Reid in \cite{wg} to construct some 3-folds and surfaces with quotient singularities  in codimension $3$ and 5. They recovered the list of 69 terminal \(\QQ\)-Fano 3-folds of Alt\i nok as a quasilinear section of \(w\Gr(2,5)\) or of a projective cone over it.   The   terminal \(\QQ\)-Fano 3-folds form a bounded family of varieties  and      more  importantly  lie in the Mori category of varieties: they are minimal models in dimension 3. Their existence in codimension 1-4 has been established in the past, for instance see \cite{ABR,tom-jerry-1}.  
Our algorithmic approach allows one to search for  terminal \(\QQ\)-Fano 3-folds in higher codimensions.  The previous attempts    \cite{qs,qs2} of constructing them as quasilinear sections of some weighted flag varieties  were not successful but  using our algorithmic approach, one can at least confirm the non-existence of  these varieties in the corresponding weighted flag varieties. Since every quotient singularity is log-terminal \cite{kawamata}, the \(\QQ\)-Fano 3-folds computed in section 4 are log-terminal. The class of log-terminal \(\QQ\)-Fano 3-folds  also forms a bounded family of algebraic  varieties \cite{borisov}. We construct some families of  isolated log-terminal  \(\QQ\)-Fano 3-folds as weighted complete intersection in the codimension 8 weighted \(G_{2}\) variety. 
\subsubsection*{Acknowledgements}
  I am grateful to Gavin Brown for some stimulating discussions which made me think about this project. I am also thankful to   Sohail Iqbal, Bal\'azs Szendr\H oi and      Shengtian Zhou for  fruitful  discussions and   Miles Reid for making  available his \textsc{Magma} function {\tt Qorb}.  Thanks are also due to Klaus Altmann and  Berlin Mathematical school for their  hospitality during  a part of this project at the  Freie Universit\"at Berlin. Last but not least, I wish to thank an anonymous referee whose comments on an earlier version helped me a great deal to improve the presentation of the article. This research is supported by the LUMS's\  faculty start up research grant.

\section{Definitions and tools }
\subsection{Baskets and polarized orbifolds}
Let the multiplicative  group  of \(r\)-th roots of unity \(\mu_r\) acts on \(\AAA^n\) by  diagonal representation \[\mu_r \ni\ep:(x_1,\ldots,x_n)\mapsto (\ep^{a_1}x_1,\ldots,\ep^{a_n}x_n).  \] Then the quotient \(\pi:\AAA^n\to \AAA/\mu_r\) is called a cyclic quotient singularity of type  \(\frac1r(a_1,\ldots,a_n)\). The cyclic quotient singularity is called isolated if \(\gcd(r,a_i)=1 \text{ for all } 1\le i\le n.\) We call a collection  (with possible repeats) of singularities \(\frac1{r_j}(a_{1_j},\ldots,a_{n_j}) \)   the  basket   of singularities, denoted by \(\sB\).   

A polarized \(n\)-fold is a pair \((X,D)\), where \(X\) is an  \(n\)-dimensional projective  algebraic variety and \(D\) is a \(\QQ\)-ample Weil divisor on \(X\). All our polarized \(n\)-folds  are well-formed and quasi-smooth; appearing as  projective subvarieties of some weighted projective space denoted by \(\PP[w_0,\ldots,w_N]\) or \(w\PP^N  \) or \(\PP[w_i]\).   We call a weighted projective variety  \(X\subset \PP[w_i]  \) of codimension \(e\)  well-formed if the singular locus of \(X\) does not contain the codimension \(e+1\) singular strata   of \(\PP[w_i].\) The subvariety \(X\subset \PP[w_i]\) is quasi-smooth if the affine cone 
\(\widetilde X =\Spec R(X,D) \subset \AAA^{N+1}\) of \(X\) is smooth outside its vertex \(\uline{0}\). 
Thus all our polarized \(n\)-folds \((X,D)\) are orbifolds, i.e. they only have quotient singularities induced by the  
singular strata of \(\PP[w_i]\). In particular, we are interested in    orbifolds  with only isolated cyclic quotient singularities. A well-formed \(n\)-fold is called projectively Gorenstein if \begin{enumerate}
\item \(H^i(X,\Oh_{X}(m))=0 \) for all \(m \text{ and }0 <i< n;\)  

\item the orbifold canonical sheaf of \(X\) is given by 
\(\om_{X}\sim\Oh_{X}(k).\) \end{enumerate} 
The integer $k$ is known as the canonical weight of  the \(n\)-fold \(X\). A variety 
\(X\) is said to have  terminal (log-terminal) singularities if, given a resolution of singularities \(f:Y \to X\), with  \[K_Y=f^*(K_X)+\sum a_iE_i,\;\; E_i \text{ are exceptional divisors and }a_i \in \QQ,\]  we have \(a_i>0\;(a_i>-1)\).   
\subsection{Graded rings and Hilbert Series} Given a  polarized \(n\)-fold \((X,D)\), the associated finite dimensional vector spaces $H^0
(X, \Oh_X(mD))$
fit  together to give rise to a finitely generated graded ring
$$ R(X, D) = \displaystyle\bigoplus_{m>0}
H^0
(X, \Oh_X(mD)),$$ with
 $$X \cong \Proj R(X, D)\into \PP[w_0,\cdots,w_N].$$ The Hilbert series of a polarized projective variety \((X,D)\), which is the Hilbert series of the graded ring \(R(X,D)\), is given by    
\[ P_{(X,D)}(t)=\sum_{m \geq 0}h^0(X,mD)\ t^{m},
\]where \(h^0(X,mD)=\dim H^0(X,\Oh_X(mD)).\)
We usually write \(P_X(t)\) for the Hilbert series for the sake of brevity. By standard 
Hilbert--Serre theorem \cite[Theorem 11.1]{atiyah}, 
$P_X(t)$ has a compact form
\begin{equation}
P_X(t)=\dfrac{H(t)}{ \displaystyle\prod_{i=0}^N\left(1-t^{w_i}\right)}.
\label{reducedhs}
\end{equation}
  The Hilbert numerator \(H(t)\) is a Gorenstein symmetric polynomial of degree \(q\): \(t^qH(1/t)=(-1)^e H(t)\) where \(e\) is the codimension of \(X\). The polynomial  \(H(t)\)  has the  form \( 1-\sum t^{b_{0j}}+\sum t^{b_{1j}}-\cdots +(-1)^c t^q\), where \(b_{0j}\) are the degrees of the equations, \(b_{1j} \) the degrees of the first syzygies, and so on. The degree \(q\) of \(H(t)\) is called the adjunction number of \(X\).   
\subsection{Weighted flag varieties}\label{sec:wfv}
Let $G$ be a reductive Lie group, with fixed Borel and maximal torus
$T\subset B \subset G$. Let $\Lam_W = \Hom(T,\CC^*)$ be the weight lattice of $G$ and let $V_\la$ denote the $G$-representation with highest
weight \(\la.\) Then there is an embedding \(\Si\into \PP V_\la\)   of a flag  variety $\Si = G/P$, where
$P = P_\la$ is the parabolic subgroup of \(G\) corresponding to the set of simple roots of $G$ orthogonal to the weight vector \(\la.\)  

Let $\Lam_W^* = \Hom(\CC^*,T)$ be the lattice of one-parameter subgroups of \(G\).  Choose \(\mu \in \Lam_W^* \) and a non-negative integer \(u \in\ZZ\) such 
that \(
\left<w\lam,\mu\right>+u >0
\label {weights} 
\)
for all  elements \(w\) of the Weyl group \(W\) of the Lie group \(G\), where $\langle ,\rangle$ denotes the perfect 
pairing between $\Lam_W$ and $\Lam_W^*$. Then we define the weighted flag 
variety \(w\Si\subset w\PP V_\la \) following  Corti--Reid \cite{wg}: take the affine
cone $\widetilde{\Sigma} \subset  V_{\lambda}$ of the embedding \(\Si \subset \PP V_\la\)  and quotient out by the $\CC^*$-action on $V_\la\backslash\{\uline 0\}$ defined by \[ 
(\varepsilon \in \CC^*) \mapsto ( v \mapsto \varepsilon^u(\mu(\varepsilon)\circ v)). 
\]  The notation \(w\Si\) will refer to general weighted flag variety and \(w\Si(\mu,u) \)  for the weighted flag variety with chosen fixed parameters \(\mu\) and \(u. \) We use the term format for each \(w\Si(\mu,u)\)  following \cite{formats}.

 \begin{theorem}\cite[Thm. 3.1]{qs} The Hilbert series of the weighted flag variety 
$\left(w\Sigma(\mu,u),D\right)$ has the following closed form.
\begin{equation}
P_{w\Si}(t)=\dfrac{\displaystyle\sum_{w\in W}(-1)^w \dfrac{t^{\left<w\rho, \mu\right>}}{(1- t^{\left<w\lambda,\mu\right>+u})}}{\displaystyle\sum_{w\in W}(-1)^w t^{\left<w\rho, \mu\right>}}.
\label{whhs}
\end{equation}
Here \(\rho\) is the Weyl vector, half the sum of the positive roots of \(G\), and  $(-1)^w=1  \; \mbox{or} -1$ depending on whether $w$ consists of an even or odd number of simple reflections in the Weyl group $W$; and \(D\) is obtained as the pullback of \(\Oh_{w\Si}(1)\) under the embedding \(w\Si \subset w\PP V_{\lam}\). 
\end{theorem}
\begin{rmk} The Hilbert series expression \eqref{whhs} always  reduces to the standard expression  \eqref{reducedhs} after performing some simplifications, see \cite{qs,qs-ahep} for details.
\end{rmk}

 \subsection{Hilbert series of isolated orbifolds}     The  integral part of our  algorithm is based on the following theorem of Buckley--Reid--Zhou \cite{brz}. The theorem gives the decomposition  of the Hilbert series of an isolated orbifold \((X,D)\) as a sum of two expressions: the initial term \(P_{I}(t)\) which   represents the smooth part of \(P_{X}(t)\)  and  the orbifold term $\displaystyle\sum_{Q_i \in\sB} P_{Q_i}(t)$ which represents the singular  part of \(P_X(t)\). 
 
\begin{theorem}\cite{brz}  \label{hs:brz}
Let $\bigl(X,D\bigr)$ be a  projectively Gorenstein quasi-smooth orbifold of dimension \(n\geq 2  \) and canonical weight \(k\), with only isolated orbifold points, given by  
$$\sB=\bigl\{Q_i \text{ of type } \frac {1}{r_i}(a_{1},\dots,a_{n})\bigr\}$$ as its only singularities. Then the Hilbert series $P_X(t) =\displaystyle\sum_{m\ge0} h^0(X,\Oh_X(mD))t^m$ of
$X$ has the form
\begin{equation}\label{ice:eq}
P_X(t)=P_I(t)+\sum_{Q_i \in\sB} P_{Q_i}(t),
\end{equation}
with { initial term} $P_I(t)$ and { each orbifold term}
$P_{Q_i}(t)$ characterised as follows:

\begin{itemize}
\item The initial term \begin{equation}\label{init-cont}P_I(t) = \dfrac{A(t)}{(1-t)^{n+1}}\end{equation} has  numerator $A(t),$ an integral Gorenstein symmetric polynomial, of
degree equal to the  coindex $c=k+n+1$ of $X$, so that $P_I(t)$ equals
$P_{X}(t)$ up to degree $\lfloor\frac c2\rfloor$ and  \(P_I=0\) for \(c < 0.\)

\item Each orbifold point \(Q_i \in \sB\)  contributes \begin{equation}\label{orb-cont}P_{Q_i}(t)=
\dfrac{B(t)}{(1-t)^n(1-t^r)}\end{equation} to the Hilbert series  where the   
numerator $$B(t)=\hbox{\rm{InvMod}}\left( \prod_{i=1}^n\frac{1-t^{a_i}}{1-t},\dfrac{1-t^r}{1-t},\left\lfloor\frac c2\right\rfloor+1     \right).$$ This is the unique integral Laurent poly\-nomial supported in
$\Bigl[ \lfloor\frac c2\rfloor+1,\lfloor\frac c2\rfloor+r-1 \Bigr]$     equal to the inverse of  \(\displaystyle\prod_{i=1}^n\frac{1-t^{a_i}}{1-t}\)  modulo \(\dfrac{1-t^r}{1-t}.\) 
 The numerator \(B(t)\) is a Gorenstein symmetric polynomial of degree \(k+n+r.\) \end{itemize}
\end{theorem}
\begin{example} Consider the Hilbert series of the  canonical 3-fold hypersurface of degree 7: $$X_7\into\PP[1,1,1,1,2].$$ $X_7$ contains an   isolated singular point of type $\sB=\left\{\frac12(1,1,1)\right\}$. The Hilbert series of $X_7$ is given by $$P_X(t)=\dfrac{1-t^7}{(1-t)^4(1-t^2)}.$$ The coindex of $X_7$ is given by  $c=1+1+3=5$. The Hilbert series  has a decomposition into the smooth part $P_I(t)$ and the orbifold part $P_Q(t)$.
The smooth part is given by $$P_I(t)=\dfrac{1+t^2+t^3+t^5}{(1-t)^4}$$ and the orbifold part is given by $$P_Q(t)=\dfrac{-t^3}{(1-t)^3(1-t^2)}.$$ The numerator of $P_Q(t)$ is $-t^3$ and we have $-t^3\equiv 1\mod (1+t)$. Thus   
$$P_X(t)=\dfrac{1+t^2+t^3+t^5}{(1-t)^4}+\dfrac{-t^3}{(1-t)^3(1-t^2)}=\dfrac{1-t^7}{(1-t)^4(1-t^2)}.$$
\end{example} 
\begin{rmk} The Hilbert series of an isolated orbifold \(X\), appearing as a   weighted complete intersection in some format \(w\Si(\mu,u)\) or cone over  \(w\Si(\mu,u)\), can be computed from the Hilbert series of the  \(w\Si\). In fact, the Hilbert numerator \(H(t)\) does not change and the denominator corresponds to the  weights of the weighted projective space containing \(X\).
\end{rmk}
\section{Algorithm to compute isolated orbifolds}
   Let   \(w\Si(\mu,u)\into \PP[w_i]\) be the  format of dimension \(d\),   corresponding to  the  fixed parameters \(\mu\) and \(u\). We aim to compute all possible candidate families of  \(n\)-folds with isolated orbifold points, whose general member is  a weighted complete intersection  of \(w\Si(\mu,u)\) or  of  projective cone(s) over it, denoted by \(\Cc^aw\Si\), having fixed canonical weight \(k\). We present an algorithmic approach  to compute  lists of such  \(n\)-folds with isolated quotient singularities   as a weighted complete intersection in \(w \Si (\mu,u)\) such that \(\Oh(K_X)=\Oh(kD)\) for \( n\ge 2 \), where \(D\) is an ample \(\QQ\)-Cartier divisor. 
 
 \subsection{ The algorithm}
 \begin{description}
\item[Step 1: Compute Hilbert series and canonical class of \(w\Si\).]We start with  a  fixed weighted  flag variety \(w\Si(\mu,u)\) and compute its  Hilbert series  \(P_{w\Si}(t)\). Each  choice of input parameters \(\mu\) and \(u\) leads to  a   codimension \(e\) embedding    \[w\Si^d(\mu,u)\into \PP^{m}[w_0,\ldots,w_m],\] where $e=m-d$. We choose the parameters $ \mu \text{ and } u$ such that \(w_i>0\), for all \(i=0,\ldots, m\). The Hilbert series of \(w\Si\) has the compact form  \begin{equation}\label{hs-wsigma}P_{w\Si}(t)=\dfrac{H(t)}{\displaystyle\prod_{i=0}^m\left(1-t^{w_i}\right) }.\end{equation}  If \(w\Si\) is well-formed, the canonical 
divisor class \(K_{w\Si}\) of \(w\Si\) is given by, \[\om _{w\Si}=\Oh_{w\Si}\left(q-\sum_{i=0}^m w_i\right)=\Oh_{w\Si}(-p).\] Since  the flag varieties and therefore weighted flag varieties are Fano varieties,    \(w\Si\) is an anti-canonically polarized variety:  \(\om_{w\Si}=\Oh_{w\Si}(-p)\),  where \(p\) is   an  integer: a multiple of \(u\) if the Lie group corresponding to \(w\Si\) is simple.  
\item[Step 2: Find all possible embeddings of \(n\)-folds \(X\) with \(\om_X=\Oh(k)\).  ] Given the embedding \[w\Si^d \into \PP[w_0,\ldots,w_m],\] we find all possible \(n\)-folds \(X\) as weighted complete intersections  inside \(w\Si \) such that \(\om_X=\Oh_X(k)\). We intersect \(w\Si\) with   generic hypersurfaces of degrees equal to some  weights $w_i$ of the weighted projective space \(\PP[w_0,\ldots,w_m] \):    
  \[X=w\Si\cap \left(w_{i_1}\right)\cap \cdots \cap\left( w_{i_l}\right)\into \PP[w_{i_0},\cdots,  w_{i_s}]\] where \(s=m-l\). We choose  \(l\) quasilinear forms  \(\left(f_j\right)\) of degree \(w_{i_j},\) such that
\(\dim(X)=d-l=n\) and $\displaystyle -p+\sum_{j=1}^l w_{i_j}=k$. 

  More generally  we can take complete intersections inside projective cones
over $ \Cc^a w\Si;$  we can   add  some more variables of degree \(w_{i}\)  to the coordinate ring which are not involved in any 
defining relation of \(w\Si\) and construct \(X\)  as its quasilinear section with \(\om_X=\Oh_X(k)\). These  newly added variables will be involved in the defining equations of \(X\) when we replace some of the variables of \(\PP^m[w_i]\) with the homogeneous forms \((f_j)\). The  Hilbert numerator $H(t)$ of the Hilbert series of $\Cc^a w\Si$ is the same as that of $w\Si$ but  we need to  multiply the denominator by $\displaystyle\prod_{_i=1}^a (1-t^{w_{i}})$, where $a$ is the number of projective cones   \(w\Si(\mu,u)\). Thus taking each  projective cone of degree \(w_{i}\)  adds \(-w_{i}\)  to the canonical class \(K_{w\Si}\). This process gives  more choices of  taking quasilinear sections of an appropriate degree.
\item [Step 3: Compute Hilbert series and the initial term  of the \(n\)-fold \(X\).]We choose one of the  \(n\)-folds  \(X\into \PP^s[w_{i_0},\ldots,w_{i_s}]\) from the list computed at Step 2, to  calculate  its Hilbert series \(P_X(t)\) and the corresponding initial term \(P_I(t)\) as given by \eqref{ice:eq}. In fact the Hilbert numerator of     the Hilbert series  \(P_{X}(t)\) of  \(X\) will be the same as in \eqref{hs-wsigma}
and the denominator changes to \(\displaystyle\prod_{j=0}^s\left(1-t^{w_{i_j}}\right).\) The initial term \(P_I(t)\), given by \eqref{init-cont}, of the Hilbert series can be computed from \(P_X(t)\) by using first \(\lfloor\frac c2\rfloor+1\) terms of \(P_X(t)\), where \(c=k+n+1\) is  the coindex of \(X\).   
\item[Step 4: Compute the isolated orbifold loci of \(w\PP^s \).] As we are only interested in computing quasi-smooth candidate orbifolds, the singularities of \(X\) are induced by the weights of the weighted projective space \(\PP^s[w_{i_j}]\).   Therefore we find  all possible contributions to the basket \(\sB \) of    \(n\)-fold isolated  cyclic quotient singularities coming from the weights of \(\PP^s[w_{i_{j}}]\). Each quotient singularity will be of type  \(\frac{1}{r_{i}}(a_{i_1},\ldots,a_{i_n})\) such that \[\gcd(r_i,a_{i_j})=1 \text{ and } \sum_{i=1}^na_{i_j}+k\equiv0 \mod r_i.\]
  Each  \(n\)-tuple of integers \((a_{i_1},\ldots,a_{i_n})\) is a sublist of the \(s\)-tuple  \((w_{i_0},\ldots,w_{i_s})\) and further we require  \(r_i>1\text{ and }r_i\in S \cup G\):  \(S=\left\{w_{i_0},\ldots,w_{i_s}\right\}\) and \(G=\left\{\gcd(S^{'}):S^{'} \subset S \right\}\).     
\item[Step 5: Compute the contributions of orbifold points to the  \(P_X(t)\).] For each isolated  orbifold point \(Q_i:=\left(\frac{1}{r_{i}}(a{i_1},\ldots,a_{i_n})\right)\) of the basket \(\sB \), we compute   its contribution \(P_{Q_i}(t)\) to the Hilbert series \(P_{X}(t)\); given by the equation  \eqref{orb-cont}. We compute the list of all possible baskets coming from the weights of \(w\PP^{s-1}\), which may contribute the Hilbert series of \(X\).  The formula  \eqref{ice:eq} can be written as      \begin{equation}\label{eq:reduced}P_X(t)-P_I(t)=\sum_{Q_i\in\mathcal B}m_iP_{Q_i}(t); \end{equation}Here \(m_i\) represents the multiplicities of the  isolated orbifold points \(Q_i\),  which remains the only unknown  in the equation \eqref{eq:reduced}. It is clear from the equation \eqref{eq:reduced} that  we  need to solve  a linear algebra  problem \begin{equation}\label{sol}[\uline P(t)]=[\uline{M}][P_{Q_i}(t)]_{ij},\end{equation}where \(P(t)=P_{X}(t)-P_I(t).\) 
\item[Step 6: Examine the  candidate \(n\)-fold \(X\). ] The last step is to compute   the coefficients  \(m_i\) appearing in \eqref{eq:reduced}. If \(m_i\geq 0\) and \(m_i\in \ZZ\) for all \(i\), then \((X,\Oh(k))\) is a suitable   candidate for a projectively Gorenstein  \(n\)-fold with isolated quotient singularities and    we  restart from step 3 by picking another model of the \(n\)-fold computed at step 2.      The actual basket of singularities of the candidate variety  \(X\) will consists  of  \(Q_i\in \sB\) such that \( m_i>0.\)
\item [Repeat: Step 3--6] We repeat steps $3-6$ for  all possible  \(n\)-fold embeddings with \(K_X=\Oh_X(k)\), computed at Step 2.

 \end{description}   
\rmk  The given algorithm essentially describes the process of finding a list of all possible families of isolated orbifolds with fixed canonical weight \(k\) whose general member may be realized as a weighted complete intersection in some prescribed format \(w\Si(\mu,u)\). But the   search for candidate varieties  inside the  weighted flag variety      \(w\Si \)  is  an infinite search, as  there is no bound  on the values of input parameters \(\mu\) and \(u\). In principle,  the algorithm does not give the complete classification of certain type of isolated orbifolds in a given weighted flag variety  but  essentially  computes the  complete list up to the certain value of the adjunction number  of the Hilbert series, which depend on the values of input parameters \(\mu\) and \(u\).  We search for the candidate varieties until  we stop getting new examples or the computer search becomes unreasonably slow due to higher degree weights.

\subsection{Implementation of the algorithm }\label{rmk-basket}
  The following remarks  describe the implementation of our algorithm in detail. 
  
\begin{enumerate} 
\item \label{three}In practice, we search for the candidate orbifolds in some chosen weighted flag variety \(w\Si\) by running a code for different    values of input parameters \(\mu\) and \(u\) . The different values of parameters \(\mu\) and \(u\) may lead to the same  set of weights on the embedding  \(w\Si\into w\PP V_\la\). For instance the two choices of  parameters, \( \mu_1=(1,-1)\) and  \( \mu_2=(-1,2)\) with \(u=3,\) give the same embedding of the weighted \(G_{2}\) variety in  \(\PP^{13}[1,2^4,3^4,4^4,5]\), see Section \ref{applications}. 

To avoid  repetition  in the computer search,  we perform the computer search with predetermined lists of  input parameters    \(\mu\) and \(u\), so that only the distinct embeddings are searched. This allows  us to check for candidate orbifolds in   distinct embeddings of \(w\Si\) in some weighted projective space. One can use a simple minded program to compute such lists of input parameters in any computer algebra system. This can essentially be  stated as the Step 0 of the  search process. 
   
\item  The choice of weights on \(w\Si\), consequently on \(X\), is determined by the values of the  parameters \(\mu\) and \(u\). We arrange the inputs and run the  code, in the order of increasing the sum of the weights in  \(\PP[w_0,\ldots,w_m]\), i.e. the sum of the weights  on the \((i+1)\)-st input  is  greater or equal than the sum on the  \(i\)-th input. This is equivalent to ordering the inputs in order of increasing the adjunction number (the degree of the Hilbert numerator \(H(t)\)) of  \(w\Si\).  
 If \(w\Si\) corresponds to a simple Lie group \(G\) then a bound on \(u\) automatically bounds the parameter \(\mu\).
\item At  Step 2 of the algorithm we compute  all possible \(n\)-folds \(X\) as weighted complete intersections  inside \(w\Si \) or of projective cone(s) over \(w\Si\) such that \(K_X=\Oh_X(k)\). Since the adjunction number \(q\) remains unchanged through the process of taking cone(s) and quasilinear sections, there are only  finitely many choices of weights \(w_{i}\) for the embedding \[X\into \PP[w_0,\ldots,w_s] \text { such that } q-\displaystyle\sum_{i=0}^sw_i=k.\] This  essentially makes the process of taking projective cones to be a  finite one, which can easily be handled by  using simple algorithms. 

The geometry behind the construction of orbifolds imposes further conditions on the choices of degree \(w_i\) of projective cones and quasilinear sections; which makes  the implementation of the algorithm    strikingly faster. The degree of the projective cones is bounded by \(w_{\max}-1\). As if the degree of a cone is greater or equal than the maximum weight \(w_{\max}\) of ambient space containing \(w\Si,\)  the newly introduced variable will not   appear in any of the defining equations of \(w\Si\) and its weighted complete intersections \(X\).  Thus a cone of  degree greater or equal  than   \(w_{\max}\) will not contribute to the orbifold part of \(P_X(t)\);  bounding the degree of the projective cone by \(w_{\max}-1\).   

Further, the degree of the forms intersecting with \(w\Si(\mu,u)\subset \PP^m[w_i]\) must be equal to one of the ambient weights of the original space containing \(w\Si(\mu,u). \) Otherwise the  process of projective cones will become redundant, in the frame work our construction. The number of projective cones  must always be less than or equal to \(s\); the dimension of  ambient weighted projective space containing \(X\).

\item   Let \[\sB=\left\lbrace m_i\times\dfrac1{r_i}\left( {a_{i{_1}},a_{i_2}},\ldots,a_{i_{n}}\right): 1\leq i \leq j\right\rbrace,\] be the  the basket of isolated orbifold points of \(X\into w\PP^s \); where \(m_i\) represents  the multiplicity of the singular point of type \(\dfrac1{r_i}\left( {a_{i{_1}},a_{i_2},\ldots,a_{i_{n}}}\right)\). Then we  define the integer \(j\) to be the  length  of the basket \(\sB\). It is evident that   \(1\le j \le b\le s\), where \(b\) is the number of non-trivial weights of \(w\PP^s\).

 Step 4 of the algorithm computes  all  possible \(n\)-fold distinct isolated orbifold  points of the weighted projective \(\PP[w_0,\ldots,w_s]\) which may potentially lie on the  candidate orbifold \(X\),   contributing to the orbifold part of the \(P_X(t)\). The total number \(B\) \(\) of such orbifold points is  usually  larger than the actual admissible length \(b\) of the basket \(\sB\) on \(X\), except when the weights of \(w\PP^s\) are relatively small.  In implementing  the algorithm, we find the set \(\mathbb{B}\) of all the admissible baskets  on \(X\); the length of the baskets  range from length 1 up to minimum\((b,B)\).  We search for the candidate orbifolds by running the code through  all the elements of \(\mathbb B\): solving  equation \eqref{eq:reduced} for each basket \(\sB\) in \(\mathbb B\).    
\item Step 6 of the algorithm checks for the solutions \(m_i\) to the equation \eqref{eq:reduced} for the given basket \(\sB\) of orbifold points. In certain cases --there may be a kernel of the singular strata--  there may exist a collection of  orbifold points \(\lbrace Q_{i}\rbrace\) in \(\sB\) such that \(\sum P_{Q_i}(t)=0\). For example, if \(Q_1=\frac15(3,3,4)\) and \(Q_2=\frac15(1,2,2)\) are two orbifold points with \(k=0\) then the corresponding orbifold terms \(P_{Q_1}(t)\text{ and }P_{Q_2}(t)\) are linearly dependent:   \[P_{Q_1}(t)= \dfrac{t^3-t^4+t^5}{(1-t)^3(1-t^5)}=-P_{Q_2}(t).\]  In the implementation of the algorithm,   we  calculate all  possible such combinations of the singular strata of the ambient weighted projective space \(w\PP^s\) containing  the candidate orbifold \(X\). Though from the routine computer search, we can only conclude that the candidate orbifold \(X\) either contains some combination of those points with each point having equal multiplicity or contain none of them. But one can precisely answer this question by explicitly computing the orbifold loci of \(X\) by using the  equations. 
\item The equation \eqref{sol} leads to some matrix equation once it is evaluated  on the the appropriate \(j\) integers, where \(j\) is the length of the corresponding basket.  The solution of the corresponding  matrix represents the multiplicities of the orbifold points of \(X\). For a given isolated orbifold \((X,D)  \) the polynomial equation of type \eqref{eq:reduced}, obtained after appropriately clearing the denominators appearing in \(P_{I}(t)\) and \(P_{Q_i}(t)\)'s,  holds for every  value of \(t\). Therefore  the corresponding matrix equation \eqref{sol} will always have a unique solution.

We are not given with an isolated orbifold; instead we   we are rather searching for the one by using the numerics of the Hilbert series and the ambient weighted flag variety. But if the numerics of the singularities and Hilbert series  correspond to some isolated orbifold in the given format then the    linear system will have a  unique solution: which are exactly the cases of our interest.  Thus the list we obtain will simply be the over list of the actual number of orbifolds in the given \(w\Si(\mu,u)\).

 \end{enumerate}

\section{Applications in 3-fold case }\label{applications}
In this section, we apply the algorithm by using a   \textsc{Magma} code, given in appendix \ref{sec:code},  to compute   lists of candidate 3-folds  with \(k\) equal to \(-1,0,\text{ and }1 \) in two types of weighted flag varieties  having embedding in codimension 3 and 8.  More precisely, we compute lists of    Fano  3-folds with isolated  log-terminal quotient singularities, Calabi--Yau 3-folds with isolated canonical quotient singularities  and canonical 3-folds with terminal quotient singularities.  We explicitly construct 5 new families of of log-terminal \(\QQ\)-Fano 3-folds as weighted completed intersection of  codimension 8 weighted \(G_2\) variety. 
\subsection{Weighted \(G_2\) variety: codimension 8}  We briefly review the construction of the codimension eight weighted flag variety;  a weighted homogeneous variety for the simple Lie group \(G_2\).  The more detailed treatment can be found in \cite{qs}.
 \begin{figure}[H]
\centering
\scalebox {1}{\input{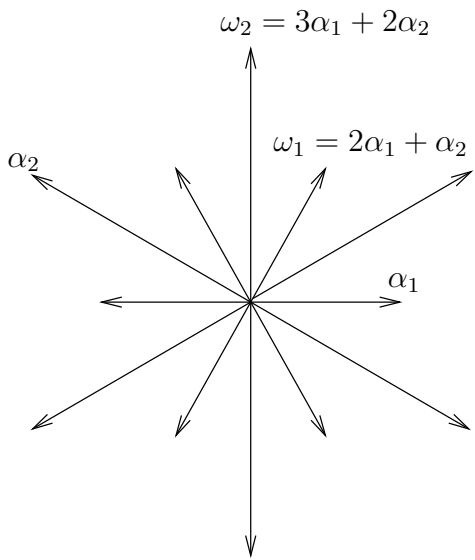}}  
\caption{Root System of \(G_{2}\)}
\label{rootsg2} 
\end{figure}
Let $\al_1, \al_2 \in \Lam_W$ 
be the  pair of simple roots of the root system $\nabla$ of the simple Lie group $G_2$.
We take  $\al_1$ to be the short simple root and $\al_2$ the long one; see Figure \ref{rootsg2}.
The fundamental weights are $ \omega_1=2\al_1+\al_2 $ and $ \omega_2=3\al_1+2\al_2$. 
The  sum of the fundamental weights; also known as Weyl vector \(\rho\) is given $\rho =5\al_1+3\al_2$.  
The cone of dominant weights is spanned by \(\om _{1}\) and \(\om_2\). Then the  $G_2$-representation with highest weight  $\lambda=\om_2=3\al_1+2\al_2$ is 14 dimensional. The corresponding  homogeneous variety 
$\Sigma \subset \PP V_{\lam} $ is five dimensional, so we have a of codimension 8 embedding \(\Si^5 \into \PP^{13}\) . Let ${\{\be_1,\be_2\}}$ be the basis of  the lattice $ \Lam_W^*$ 
; dual to \(\{\al_1 , \al_2\}\).

The weighted version can be constructed by taking $\mu=a\be_1+b\be_2 \in \Lam_W^*$ and \(u\in \ZZ^+\); following  Section \ref{sec:wfv}. The defining equations of the  weighted \(G_2\) variety \(w\Si\) can be calculated by using the decomposition of the second symmetric power \(S^2(V_\la^*)\) of the dual representation of \(V_\la\), as a module over the Lie algebra \(\gog_2\): see \cite{qs} for the details. The   28 quadrics cut out the defining  locus of \(w\Si\), explicitly given in \cite[Appendix A]{qs}. 

   The weights of the ambient weighted projective space \(w\PP^{13}\) are parameterized by an integer value vector \(\mu=(a,b)\) and \(u\);  given by  \begin{equation}\label{g2:wts}[\pm a+u,\pm b+u,\pm (a+b)+u,\pm (2a+b)+u,\pm (3a+b)+u,\pm (3a+2b)+u,u,u].\end{equation}  As we need all the weights to be positive; there are   finite choices of the parameters  \(\mu\) for each positive integer \(u \). The  adjunction number \(q\), the degree of the Hilbert numerator \(H(t)\), is \(11u\). The sum of the weights is \(14u\); if  \(w\Si\) is well-formed then the canonical divisor class is \(K_{w\Si}=\Oh_{w\Si}(-3u)\). We search for the examples in order of increasing the sum of the weights \(\sum w_i=14u\) on \(w\PP V_\la\), which  corresponds to the increase in  \(u\) and consequently in order of increasing the adjunction number \(q=11u\).

\subsection{Examples and lists of isolated orbifolds in codimension 8}
In this section we prove the existence of some families of isolated log-terminal \(\QQ\)-Fano 3-folds whose general member can be embedded in a codimension 8 weighted \(G_2\) variety. We also present the list of examples obtained for the Calabi--Yau 3-folds and canonical 3-folds in the codimension 8 weighted \(G_2\) variety.   

\thm \label{thrm:g2fano}Let \(w\Si\) be the codimension eight  weighted \(G_2\)-variety.  Then there exist 6 families of isolated log-terminal \(\QQ\)-Fano 3-folds whose general member is a weighted complete intersection in \(w\Si\) or  some projective cone(s) over \(w\Si\), given by the Table    \ref{tab:fano3g2}.
\begin{table}[H]

\begin{center}\caption{ Log-terminal $\QQ$-Fano 3-folds in weighted $G_2$ variety, \((\mu,u) \)-Input parameters, Weights-- weights of weighted projective space containing \(X\), \((-K_X)^3\)--Degree of \(X\),  Basket-- orbifold points of \(X\), and BK-- Basket kernel for \(X\). }\label{tab:fano3g2}
\renewcommand{\arraystretch}{1.5}
\begin{tabular}{|c| c| c| c|c| }  

\hline\hline 
  \((\mu:u)\)  & Weights   & \((-K_X)^3\)& Basket&BK\\\hline (0,0:1) &\(\PP\left[1^{12}\right]\)& \(18\)&\( \)&N\\ 
\hline (-1,1:3) &\(\PP\left[1,2^4,3^4,4^2,5\right]\)& \(\frac{9}{10}\)&\(9\times\frac12(1,1,1),\frac 15(3,4,4) \)&Y\\
\hline (-1,1:4) &\(\PP\left[2,3^{4},4^3,5^3\right]\)& \(\frac{1}{5}\)&\(2\times\frac12(1,1,1),6\times \frac13(1,1,2),3\times\frac 15(3,4,4)\)&N\\
\hline (-2,3:4) &\(\PP\left[1^2,2,3^2,4^{3},5^2,6,7\right]\)& \(\frac{9}{14}\)&\( 2\times \frac 14(1,1,3),\frac 17(4,5,6)\)&N\\
\hline (-4,6:7) &\(\PP\left[1^2,3,5^2,7^3,9^2,11,13\right]\)& \(\frac{18}{91}\)&\( 2\times \frac 17(1,2,5),\frac{1}{13}(7,9,11)\)&N\\
\hline (-3,4:7) &\(\PP\left[2,3,4,5,6^2,7^2,8,9,10,11\right]\)& \(\frac{4}{65}\)&\( 7\times \frac 12(1,1,1),3\times \frac13(1,1,2),\frac 1{11}(6,7,10)\)&Y\\
\hline 
\end{tabular}
\end{center}
\end{table}
\proof The proof is essentially to construct  all the candidate examples listed in Table \ref{tab:fano3g2}, one after another. We start with the first entry listed in Table \ref{tab:fano3g2}, then we have the following data to begin with:\label{log:G2}

\begin{itemize}
\item Input: $ \mu=(-1,1) $, $ u=3 $
\item Variety and weights: $ w\Si_{} \subset \PP^{13}[1,2^4,3^4,4^4,5]  $
\item Canonical class: $ K_{w\Si}= \Oh(-9) $, as we can check that $w\Si $ is well-formed.
\item Hilbert Numerator: $1-3 t^4-6 t^5-8 t^6+6 t^7+21 t^8+\ldots+6 t^{26}-8 t^{27}-6 t^{28}-3 t^{29}+t^{33}$
\end{itemize}  
We take a threefold quasilinear section, by intersecting  \( w\Si\) with two  general forms of degree four   \[X=  w \Si\cap(4)^2\into \PP^{11}[1,2^4,3^4,4^2,5] \text{ with } K_X=\Oh(-1).\] 
We check all the singular strata of \(X\) by using the equations given in \cite[Appendix A]{qs}.\\[1.5mm]
\( 1/5 \) \textbf{singularities}: Since we are not taking any degree five sections of \(w\Si\) and the variable of weight 5 does not appear as monomial of degree 2 in the defining equations of \(w\Si\),  \(X\) contains this point. By using the implicit function theorem we find the local transverse parameters near this point to be of degree 4,4 and 3. Therefore \(X\) has a singular point of type \(\frac 15(3,4,4)\). \\[1mm] 
\(1/4\) \textbf{singularities}: \(X\) does not contain this singular point.  \\[1mm]
\(1/3\)\textbf{ singularities}:  \(X\) also avoids  \(\)  singularities of type \(\frac13\).\\[1mm]
\(1/2\)  \textbf{singularities} \(X\) contains 9 singular points of type \(\frac 12(1,1,1).\)   \\[1mm] 
Thus \((X,-K_{X})\)  is a log-terminal  \(\QQ\)-Fano threefold with $(-K_X)^3= \frac{9}{10}$  and with \(9 \times \frac 12(1,1,1), \frac 15(3,4,4)\) singular points.
 
 The existence of  the remaining    log-terminal \(\QQ\)-Fano 3-folds appearing in Table \ref{tab:fano3g2}, has been established by using explicit equations and a similar analysis of the singular strata.\(\square\)
\rmk If \(X\) is a terminal \(\QQ\)-Fano 3-fold appearing as a quasilinear section of some weighted \(G_{2}\) variety with parameters \((\mu:u)\), then the adjunction number of \(X\)  will be \(11u.\) The list of all possible codimension 8 terminal \(\QQ\)-Fano 3-fold, with their  basket of singularities, adjunction number, degree etc; is available on the graded ring database page \cite{grdb}. The  highest adjunction number for the codimension 8 terminal \(\QQ\)-Fano 3-folds  which is an integer multiple of 11 is 66 in the database, so a terminal \(\QQ\)-Fano 3-fold in the weighted \(G_2\) variety can exist only for  \(u\le 6\).   But the   computer search   does not produce even  a candidate terminal \(\QQ\)-Fano 3-fold for \(1\le u \le6\), leading to the following corollary.
 \begin{coro}
There does not exist a terminal \(\QQ\)-Fano 3-fold in codimension 8  which can be realized  as  weighted complete intersection in codimension 8  
weighted \(G_2\) variety. \end{coro}

\rmk  Table \ref{tab:fano3g2}  does not represent a complete classification  of  isolated  log-terminal \(\QQ\)-Fano 3-folds in the codimension 8 weighted \(G_2\) variety. Certainly, it is a sublist of the  complete classification of such 3-folds in the given weighted \(G_{2}\) variety, as the computer search has been completely performed only up to \(u=7\). One should expect more examples for the higher values of \(u\).    Similarly Table \ref{tab:c3f} and Table \ref{tab:cy3} do not represent the full classification of Canonical and Calabi--Yau 3-folds respectively.
 On the other hand  the computer routine  returns an over list of  numerical candidates,  up to certain values of the adjunction number \(q\).  In the case of log-terminal \(\QQ\)-Fano 3-folds, we get a list of  33  suggested models of candidate orbifolds up to \(q=77\) but only 6 of them exist as an actually variety with the suggested invariants and singularities. We include the smooth Fano 3-fold of genus 10 in the list for the sake of completeness as it already appeared in \cite{qs}. The existence of these candidates is established by checking through the equations of each of these 33 numerical candidates. The checking is done partly by hand and partly by the computer algebra, as the number of defining equations is substantially large. 
  The following example represents a suggested candidate  orbifold which does not exists as an actual variety with the suggested numerical data. 
\begin{example}
If we take the \(\mu=(-3,4)\) and \(u=6\), then we get the embedding \[w\Si(\mu,u)\into\PP\left[1,2,3,4,5^2,6^2,7^2,8,9,10,11\right],\] with \(K_{w\Si}=\Oh(-18).\) Then a 3-fold quasilinear section \[X=w\Si\cap(10)\cap(7)\into\PP\left[1,2,3,4,5^2,6^2,7,8,9,11\right]\] has \(K_X=\Oh_X(-18+10+7)=\Oh_X(-1). \) The computer search suggests a numerical candidate   log-terminal \(\QQ\)-Fano 3-fold with the basket of singularities \(\sB=\left\{\frac{1}{11}(6,8,9),\frac15(1,1,4), 2\times \frac12(1,1,1)\right\}\) and \((-K_X)^3=\frac{9}{55}\).
By checking through the equations of \(X\), induced from the equations of \(w\Si\), the orbifold points of type \(\frac1{11}(6,8,9)\) and \(\frac15(1,1,4)\) lie on \(X\). But the fix locus under the action of cyclic group \(\ZZ_2\) is  an empty subscheme of \(X\); the 2 points of type \(\frac12(1,1,1) \) do not actually lie on \(X\). Thus the suggested candidate 3-fold model can not be realised as a quasilinear section of weighted \(G_2\) variety with the given numerical data.
\end{example}
\begin{table}[h]
\begin{center}\caption{ Candidate canonical  3-folds  in the codimension 8 weighted \(G_2\) variety  }\label{tab:c3f}
\renewcommand{\arraystretch}{1.5}
\begin{tabular}{| c| c| c| c|}  
\toprule\hline 
  \((\mu:u)\)  & Weights   & \((K_{X})^3\) & Basket\\
\hline (0,0:2) &\(\PP\left[1^3,2^9\right]\)& \(9\)&\(18\times \frac12(1,1,1)\)\\
 
 \hline \multirow{3}{*}{(-1,1:3)} &\(\PP\left[1^2,2^4,3^3,4^2,5\right]\)& \(\frac{27}{10}\)&\(9\times \frac12(1,1,1),\frac 15(1,4,4) \)\\ \cline{2-4}
&\(\PP\left[1,2^5,3^4,4^3\right]\)& $\frac95$& $18\times \frac12(1,1,1),\frac15(2,3,4)$\\
\cline{2-4}
&\(\PP\left[1,2^4,3^5,4^2\right]\)& $\frac95$& $9\times \frac12(1,1,1),6\times\frac13(1,2,2)$\\ \hline
\multirow{2}{*}{(-2,3:4)}&\(\PP\left[1^{2},2^2,3^2,4^2,5^2,6,7\right]\)& $\frac97$& $2\times \frac12(1,1,1),\frac17(2,5,6)$\\
\cline{2-4}
&\(\PP\left[1,2^2,3^3,4^2,5^3,6\right]\)& $\frac35$& $2\times \frac12(1,1,1),3\times \frac13(1,2,2),2\times\frac15(2,3,4)$\\\hline
\multirow{2}{*}{(-1,1:4)}&\(\PP\left[1,2,3^3,4^4,5^3\right]\)& $\frac35$& $2\times \frac12(1,1,1),3\times\frac15(1,4,4)$\\
\cline{2-4}
&\(\PP\left[2^2,3^4,4^3,5^3\right]\)& $\frac25$& $4\times \frac12(1,1,1),6\times \frac13(1,2,2),3\times\frac15(1,4,4)$\\
\hline
\multirow{2}{*}{(-2,3:5)}&\(\PP\left[2^2,3^2,4^2,5^3,6,7,8\right]\)& $\frac{9}{40}$& $11\times \frac12(1,1,1),2\times\frac15(1,4,4),\frac18(3,5,7)$\\
\cline{2-4}
&\(\PP\left[2,3^3,4^2,5^3,6,7^2\right]\)& $\frac{6}{35}$& $6\times \frac13(1,2,2),2\times\frac15(1,4,4),\frac17(3,4,6)$\\
\hline 
(-2,3:6)&\(\PP\left[2,3^2,4,5^2,6^2,7^2,8,9\right]\)& $\frac{1}{9}$& $2\times \frac12(1,1,1),7\times \frac13(1,2,2),\frac19(2,7,8)$\\ \hline
 (-3,4:7)&\(\PP\left[2,3,4^2,5,6,7^2,8,9,10,11\right]\)& $\frac{3}{44}$& $10\times \frac12(1,1,1), \frac14(1,3,3),\frac1{11}(4,7,10)$\\\hline
\bottomrule
\end{tabular}
\end{center}
\end{table}
 \rmk Since the   set of
orbifold contributions for the canonical terminal points
$\frac1r(-1,a,-a)$ on 3-folds with $k=1$ are linearly independent \cite{fletcher1}, we do not include the basket kernel column in Table 2. 
\begin{table}[H]
\begin{center}\caption{ Candidate  canonical Calabi--Yau 3-folds in the codimension 8 weighted \(G_2\) variety  }\label{tab:cy3}
\renewcommand{\arraystretch}{1.5}
\begin{tabular}{| c| c| c| c|c|}  
\toprule\hline 
  \((\mu:u)\)  & Weights   & \(D^3\) & Basket&Ker\\
 \hline (-1,1:3) &\(\PP\left[1,2^4,3^5,4,5\right]\)& \(\frac{6}{5}\)&\(6\times \frac13(1,1,2),\frac 15(3,3,4)\)&\multirow{5}{*}{Y}\\
\cline{1-4} (-2,3:4) &\(\PP\left[1^{2},2,3^{3},4^2,5^2,6,7\right]\)& \(\frac{6}{7}\)&\(3\times \frac13(1,1,1),\frac 17(3,5,6)\)&\\
\cline{1-4} (-2,3:5) &\(\PP\left[2^2,3^2,4,5^3,6^2,7^2\right]\)& \(\frac{6}{35}\)&\( 3\times \frac 13(2,2,2),2\times\frac 15(1,2,2),\frac 17(2,6,6)\)&\\
\cline{1-4} \multirow{2}{*}{(-3,4:7)} &\(\PP\left[2,3,4,5^2,6,7^2,8,9,10,11\right]\)& \(\frac{3}{55}\)&\( 2\times \frac 15(1,2,2),\frac{1}{11}(5,7,10)\)&\\
\cline{2-4}  &\(\PP\left[3,4,5^2,6^{2},7^3,8,9,10\right]\)& \(\frac{1}{35}\)&\( 7\times \frac 12(1,1,1),3\times \frac13(1,1,2),\frac 1{11}(6,7,10)\)&\\
\hline \end{tabular}
\end{center}
\end{table}
\subsection{Computations of  other known lists}
There are some famous lists
of 3-dimensional orbifolds, for example 95 \(\QQ\)-Fano 3-fold hypersurfaces in  weighted projective spaces \cite{fletcher} or 69 families of codimension 3 \(\QQ\)-fano 3-folds in the weighted Grassmannian \(w\Gr(2,5)\) format  etc. The summary of such lists and corresponding references can be found in Table 1 of \cite{formats}; where the results were obtained by using a  different approach than ours. Since weighted projective spaces are a particular type of weighted flag varieties; those lists of orbifolds can be recovered by  making  slight modifications to our computer routine, except  the lists of  codimension 1 Calabi--Yau 3-folds of \cite{kre}.  Theoretically, we will  eventually obtain the full lists of  these 3-folds as well    but  the computer search gets hopelessly slow once the weights of the ambient space get larger.

 In theory, the algorithm can be used to find lists of orbifolds inside any ambient  weighted projective variety with a computable canonical divisor class and Hilbert series. As a  description, we recover   lists of canonical, Calabi--Yau and log-terminal \(\QQ \)-Fano 3-folds (including the terminal \(\QQ\)-Fano 3-folds) inside weighted Grassmannian $w\Gr(2,5)$ format and present the results in table \ref{tab:summary}.
   
   Table \ref{tab:summary} presents the summary of the results obtained in two cases: codimension 8 weighted \(G_2\) variety and codimension 3  weighted Grassmannian \(\Gr(2,5)\). In each case the results are searched up to the  adjunction number \(q_{\max}\). The number  \(q_\mathrm {res}\)  represents the adjunction number for which the last numerical 3-fold was found. The column \(\#w\Si_{\rm dis}\) gives the number of distinct embeddings searched in the given format, and $\#w\Si_{\rm res}$ is the number of the embedding where  the last result was found. The column \(\#\)output represents the number of suggested candidate 3-folds and \#result gives the number of plausible candidate  3-folds in the given weighted flag variety.

In the case of \(w\Gr(2,5)\), we recover the list of 18 canonical for  3-folds computed in \cite{formats} for \(k=1\). The  list of famous 69 families of \(\QQ\)-Fano 3-folds in codimension 3 is obtained as a sublist of 403 numerical examples of log-terminal \(\QQ\)-Fano  3-folds computed for \(k=-1\).  In the case of \(k=0\) the number of Calabi--Yau 3-folds obtained up to adjunction number \(71\) is   an over list of the 187 such 3-folds appearing on the graded ring database page \cite{grdb}, computed by using a different approach than ours.

For the   codimension 8 weighted \(G_2\) variety the adjunction number increases quite rapidly;  leading to fewer cases of distinct embeddings.   For \(k=1  \)  we obtained a list of 14 numerical candidates with 12 plausible examples and the case of \(k=0\) gives 13 candidate families of Calabi--Yau 3-folds with 6 of them being plausible. For \(k=-1\) we get 33 numerical candidates of isolated log-terminal \(\QQ\)-Fano 3-folds with 6 of them existing as actual varieties;  explicitly constructed in Theorem \ref{thrm:g2fano}.       
 
 \begin{table}[t]
 
\caption{The summary of results showing the number of families of log-terminal \(\QQ\)-Fano 3-folds,  Calabi--Yau 3-folds, and canonical 3-folds with isolated orbifold points in two  formats: codimension 8 \(G_2\), and  \(\Gr(2,5)\).  The column $\#w\Si_{\rm dis}$  gives the number of distinct embeddings searched for examples, \( \#w\Si_{\rm last}\) is the last embedding where the example appeared,  $q_{\mathrm{res}}$ gives the largest adjunction number for which a result was found; $q_{\mathrm{max}}$  gives the largest adjunction number  searched;  \#outputs gives the number of candidates found by the computer;  \#results gives the number of candidates after removing the candidates with obvious failure.
\label{tab:summary}}
\[
\renewcommand{\arraystretch}{1.2}\begin{array}{|c|c|c|c|c|c|c|c|c|}
\toprule\hline
\text{Format} & \text{codim} & k & \#w\Si_{\rm dis}  &\# w\Si_{\rm{res}} & q_{\mathrm{max}}&q_{\mathrm{res}} &   \text{\#outputs} & \text{\#results} \\\hline
 \multirow{3}{*}{$G_2$} &\multirow{3}{*} {8} & -1 & 23 &19  &77   &77  &  32 & 6 \\\cline{3-9}
 &  & 0 &41  &  27&99   &  88& 12  & 6 \\\cline{3-9}
 &  & 1 &53  & 19 &110   &77  &  14 & 12 \\\hline
 \multirow{3}{*}{$\Gr(2,5)$} &\multirow{3}{*} {3} & -1 & 17180 &13403  & 63  &  63& 403(69)  &   \\\cline{3-9}
 &  & 0 & 29941 &29165  &71   & 71 & 221  & 187 \\\cline{3-9}
 &  & 1 & 29941 &1196  &  71  &35  & 18  &18  \\\hline
 
\bottomrule
\end{array}
\]
\end{table}
\rmk There are some other types of weighted flag varieties constructed in \cite{qs,qs2} in codimension \(6,7, \text{ and }9\). The detailed list of isolated orbifolds in those weighted flag varieties will appear elsewhere \cite{BKQ}. In this article,  we restrict the  attention to description of  the algorithm and its applications to sample cases.

\appendix 
\section{\textsc{Magma} code to compute $n$-folds}\label{sec:code} 
This code consists of the main \textsc{Magma} function {\tt "Format"}, which uses some auxiliary { functions} and some extra data to produce the required lists of examples.  The following is the most general form of the implementation of our algorithm. The search process can be significantly fastened by a slight modification in the code in particular cases. For example, in the case of isolated Calabi--Yau 3-folds the index of singularity must be odd, so a minor modification in the function {\tt " Porb\_Cont"} fastens the search significantly.

For the whole calculation we run the following  basic commands to start calculations after logging into  \textsc{Magma}.
\begin{verbatim}

Q:=Rationals();
R<t>:=PolynomialRing(Q);
K:=FieldOfFractions(R);
S<s>:=PowerSeriesRing(Q,50);\end{verbatim}

The function {\tt Qorb} calculates the contribution \(P_{Q_i}(t)\) of each isolated singular point \(\frac1{r_i}(a_1,\ldots,a_n)\) to the Hilbert series \(P_X(t)\) of \(X\). The input to this function are the index of singularity \(r\), the   weights of the local coordinates \({\tt LL}=[a_1,\ldots,a_n]\) and the canonical weight {\tt kx} of \(X\). This function is the
 own implementation of their algorithm by the authors of \cite{brz}.
\begin{verbatim}
function Qorb(r,LL,kx)
L := [ Integers() | i : i in LL ]; 
        if (kx + &+L) mod r ne 0
        then error "Error: Canonical weight not compatible";
        end if;
n := #LL; Pi := &*[ R | 1-t^i : i in LL];
h := Degree(GCD(1-t^r, Pi));  l := Floor((kx+n+1)/2+h);
de := Maximum(0,Ceiling(-l/r)); m := l + de*r;
A := (1-t^r) div (1-t); B := Pi div (1-t)^n;
H,al_throwaway,be:=XGCD(A,t^m*B);
return t^m*be/(H*(1-t)^n*(1-t^r)*t^(de*r));
end function;
\end{verbatim}

The function {\tt  "Init\_Term"}\    computes the contribution of the initial term \(P_I(t)\) of the Hilbert series \(P_X(t) \) of \(X\), as given by \eqref{ice:eq}. The input is: {\tt hs}- the  Hilbert series of \(X\), {\tt n}- dimension of \(X\), and {\tt c-}coindex of \(X\).  
 
\begin{verbatim}
function Init_Term(hs,c,n)
co:=Coefficients(S!hs)[1..Floor(c/2)+1];
f:=&+[co[i]*t^(i-1): i in [1..#co]];
pp:=R!(f*(1-t)^(n+1));
        if IsEven(c) eq true then
        return (&+[Coefficient(pp, i )*(t^i+t^(c-i)):i in [0..c div 2-1]]+
        Coefficient(pp,c div 2)*t^(Floor(c/2)))/(1-t)^(n+1);
        else
        return &+[Coefficient(pp,i)*(t^i+t^(c-i)):i in [0..Floor(c/2)]]
        /(1-t)^(n+1);
        end if;
end function;
\end{verbatim}

The function {\tt  "Pos\_Wt"} computes  all  possible embeddings  of \(X \into w\PP^{s-1}\), as a quasi-linear section of \(w\Si \into \PP L\) and of all possible projective  cone(s) over it, with a desired  canonical   divisor class. As an input we use the weights of the embedding  \(w\Si\into\PP[w_i]\) as a list of integers \(L\) and integers \(s,w\), where \(w\) is the required sum of the weights on \(\PP^{s-1}.\)
As an output we get  lists of integer lists of length \(s\) such that their sum is \(w\) and corresponding ambient weighted projective space  \(\PP^{s-1}[w_i]\) is well-formed.
  
\begin{verbatim}
function Pos_Wt(L,s,w)
PosWt:=[Sort(p): p in RestrictedPartitions(w,s,{1..Max(L)})| Multiplicity(p,Max(L)) 
le Multiplicity(L,Max(L)) and (&+[GCD(Remove(p,i)): i in [1..#p]] eq s)];
return PosWt;
end function;
\end{verbatim}
The function {\tt  "Porb\_Cont"} calculates   all the possible singularities coming from the  embedding of \(X\). As an input it takes the list of weights of \(w\PP^{s-1}\), the canonical weight {\tt kx} and the dimension \(n\) of \(X\).  

\begin{verbatim}
function Porb_Cont(weights,kx,n)
LL:=[PowerSequence(Integers())|]; R:=[Integers()|];
for r in weights do
gcds:={GCD(r,s) : s in weights};
weights cat:= [p : p in gcds | p ne 1 and p notin weights];
end for;
set:=SequenceToSet(Sort(weights)); 
        for r in {w : w in set | w ne 1} do
        rowt:=Sort([a mod r : a in weights | GCD(a,r) eq 1]); 
                if (#rowt ge n ) then
                N:=SetToSequence(Subsets({1..#rowt},n));
                I:=[SetToSequence(N[i]): i in [1..#N]];
                BI:=[Sort(rowt[I[j]]): j in [1..#N]|(m mod r eq 0) where 
                m is &+rowt[I[j]]+kx];
                tr:= Setseq(Seqset(BI));
                        for L in tr do
                        LL:=Append(LL,L);R:=Append(R,r);
                        end for;
                end if;
        end for;
return R,LL,weights;
end function;   \end{verbatim}

The function {\tt  "Baskets"}  computes all possible baskets which may lie on \(X\), induced from the weights of  embedding \(\PP^{s-1}[w_i]\). The  input of this function is  the output from the function {\tt "Porb\_Cont"}. we know that \(X\into \PP^{s-1}[w_i]  \), so the maximum length (defined in the Section  \ref{rmk-basket}) of  the baskets is \(s\). 
  
\begin{verbatim}
function Baskets(R,LL,weights)
RR:=[PowerSequence(Integers())|];
BB:=[PowerSequence(PowerSequence(Integers()))|];
for s in [1..Min(#R,#weights)] do
   for I in {<R[Sort(SetToSequence(m))],LL[Sort(SetToSequence(m))]> :
   m in Subsets({1..#R},s)} do
        if #(SequenceToMultiset(I[1]) meet
        SequenceToMultiset(weights)) ge s then
        Append(~RR,I[1]);
        Append(~BB,I[2]);
        end if;
   end for;
end for;
return RR,BB;
end function; 
\end{verbatim}
The function {\tt "Bask\_Kernel"} computes the kernel of the basket of singularities \(\sB \) induced from the weights of \(\PP^{s-1}[w_i]\), i.e. given a set of  isolated orbifold points \(Q_{i}\), it computes  all possible   combinations of   \(Q_{i}\)  such that \(\sum P_{Q_i}(t)=0\). As an input it takes the weights of the ambient space, the output from the function {\tt Baskets}, and the canonical weight {\tt kx} of \(X\).
\begin{verbatim}
function Basket_Kernel(weight,RR,BB,kx)
if #RR ge 2 then
  BR:=[PowerSequence(Integers())|];BL:=[PowerSequence(PowerSequence(Integers()))|];
   for s in [2..Min(#weight,#RR)] do
      for I in {<RR[Sort(SetToSequence(m))],BB[Sort(SetToSequence(m))]> :
      m in Subsets({1..#RR},s)} do
                if (#(SequenceToMultiset(I[1]) meet SequenceToMultiset(weight)) eq s) 
                and (&+[Qorb(I[1][j],I[2][j],kx): j in [1..#I[1]]] eq 0) then 
                Append(~BR,I[1]); Append(~BL,I[2]); 
                end if;
      end for;
  end for;
end if;
return BR,BL;
end function;\end{verbatim}
The procedure {\tt  "Format"} is the main function utilizing the rest of the functions in combination to search for a suitable candidate orbifold \((X,D)\).  As an input it takes the following data of \(w\Si(\mu,u)\) and \(X\).
\begin{itemize}
\item {\tt num}- The Hilbert numerator of \(w\Si(\mu,u)\) and hence of \(X\)
\item {\tt wtsigma}- The weights of the embedding \(w\Si \into w\PP^{}\), as a list of integers
\item {\tt n}- The dimension of desired candidate orbifold \((X,D)\)
\item {\tt kx}- is the canonical weights of \(X\)
\item {\tt s}- is the number coming from \(X\into w\PP^{s-1}\)   
\end{itemize}
This {\tt function} checks for  plausible candidate
isolated  orbifolds in a given weighted flag variety \(w\Si(\mu,u)\) and  all possible projective cones over it. If such a candidate is found it will output the candidate with its degree, singularities, weights of the embedding and Hilbert numerator, otherwise only the Hilbert numerator is returned to show the completion of the process.

 \begin{verbatim}
function Format(num,kx,s,wtsigma,n)
for weight in Pos_Wt(wtsigma,s,Degree(num) - kx) do
den:=&*[1 - t^n : n in weight];
px:=num / den;
deg:=Evaluate(px * (1 - t)^(n + 1),1);
pini:=Init_Term(px,n + kx + 1,n);
    if px eq pini then
    printf "Smooth %o-fold with canonical class K_X=O(%o)\n" cat
    "Ambient: P%o\nDegree: %o\n",
    n,kx,weight,deg;
    else
    Ra,LLa,weights:=Porb_Cont(weight,kx,n);
    RR,BB:=Baskets(Ra,LLa,weights);
        for i in [j : j in [1..#RR] | #RR[j] ne 0] do
        R:=RR[i];
        B:=BB[i];
        Porb:=[Qorb(R[j],B[j],kx) : j in [1..#R]];
        den:=&*[Denominator(PQ) : PQ in Porb] *Denominator(px) * Denominator(pini);
        PX:=den * px;
        Pini:=den * pini;
        Porb:=[PQ * den : PQ in Porb];
        LHS:=PX - Pini;
           if Min([Degree(LHS) - Degree(PQ) : PQ in Porb]) ge 0 then
           V:=Vector(Integers(),[Evaluate(LHS,k) : k in [2..#R + 1]]);
           A:=Matrix([PowerSequence(Integers())|[Evaluate(PQ,k):k in [2..#R+ 1]]: 
           PQ in Porb]);
           ok,sol:=IsConsistent(A,V);
                if ok  and Min(Eltseq(sol)) ge 0 and 
                (LHS - &+[sol[m]*Porb[m]:m in [1..#R]] eq 0) then
                BR,BL:=Basket_Kernel(weight,Ra,LLa,kx);
                        if  #BR ge 1 and [&+Sum,&*Prod,&+Ind] notin SS then
                        printf "Isolated %o-fold with canonical class K_X=O(%o)\n" cat
                        "Ambient: P%o\nSingularities:\n %o x 1/%o\n%o\n" cat
                        "Degree:%o\n\n\n\n" cat
                        "Basket Kernel:1/%o x %o\n\n",
                        n,kx,weight,sol,R,B,deg,BR,BL;                        
                        else 
                        printf "Isolated %o-fold with canonical class K_X=O(%o)\n" cat
                        "Ambient: P%o\nSingularities:\n %o x 1/%o\n%o\n" cat
                        "Degree:%o\n\n\n\n",
                        n,kx,weight,sol,R,B,deg;
                        end if;
                 end if;                     
            end if;
        end for;
    end if;
end for;  return "Hilbert Numerator:" num;
end function;
\end{verbatim} 
\bibliographystyle{amsplain}
\bibliography{imran}
\noindent {\sc Department of Mathematics,}

\noindent {\sc Lums School of Science and Engineering,}

\noindent {\sc U-Block, DHA, Lahore, Pakistan}
\\
\noindent {\it Email address}: {\tt i.qureshi@maths.oxon.org}

\end{document}